\documentclass[11pt]{amsart}
\usepackage{amsmath, amssymb, amsthm, fullpage, epic, eepic}
\usepackage[all]{xypic}
\newtheorem{thm}[equation]{Theorem}
\newtheorem{prop}[equation]{Proposition}

\newtheorem{cor}[equation]{Corollary}
\newtheorem{lemma}[equation]{Lemma}

\theoremstyle{definition}

\theoremstyle{remark}

\newtheorem{rem}[equation]{Remark}

\makeatletter
\renewcommand{\subsection}{\@startsection{subsection}{2}{0pt}{-3ex
plus -1ex minus -0.2ex}{-2mm plus -0pt minus
-2pt}{\normalfont\bfseries}}
\renewcommand{\subsubsection}{\@startsection{subsubsection}{2}{0pt}{-3ex
plus -1ex minus -0.2ex}{-2mm plus -0pt minus
-2pt}{\normalfont\bfseries}} \makeatother

\numberwithin{equation}{section}

\newdir^{ (}{{}*!/-5pt/@^{(}}

\newcommand{\erem}{\hfill$\lozenge$\end{rem}\vskip 3pt }

\newcommand{\iso}{{\;\stackrel{_\sim}{\to}\;}}

\newcommand{\beq}{\begin{equation}\label}
\newcommand{\eeq}{\end{equation}}

\newcommand{\onto}{\twoheadrightarrow}

\newcommand{\dep}{\operatorname{dp }}

\newcommand{\en}{\enspace }

\def\R{\mathbb{R}}

\def\Z{{\mathbb Z}}

\def\1{{1}}

\begin{document}
\title{\qquad\en On dominance and minuscule Weyl group elements}
\author{Q\"endrim R. Gashi and Travis Schedler}
\maketitle

\begin{abstract}
  Fix a Dynkin diagram and let $\lambda$ be a coweight. When does
  there exist an element $w$ of the corresponding Weyl group such that
  $w$ is $\lambda$-minuscule and $w(\lambda)$ is dominant? We answer
  this question for general Coxeter groups. We express and prove these
  results using a variant of Mozes's game of numbers.
\end{abstract}

\section{Introduction}

Mazur's Inequality \cite{mazur, mazur2} is an important $p$-adic
estimate of the number of rational points of certain varieties over
finite fields. It can be formulated in purely group-theoretic terms,
and the classical version can be viewed as a statement for the group
$GL_n$ (see \cite{kot}). Kottwitz and Rapoport formulated a converse to
this inequality \cite{krapo}, which is also related to the
non-emptiness of certain affine Deligne-Lusztig varieties, and they
reduced the proof to a purely root-theoretic problem, which is solved
in \cite{qendrim3}. A crucial step in \cite{qendrim3} involves the use
of Theorem \ref{question} below, which we state after introducing some
standard notation and terminology.

Let $\Gamma$ be a simply-laced Dynkin graph, with corresponding simple
roots $\alpha_1, \ldots, \alpha_n$, positive roots $\Delta_+$, Weyl
group $W$, and simple reflections $s_1, \ldots, s_n \in W$.  Let
$P_\Gamma$ be the lattice of coweights corresponding to $\Gamma$.
Following Peterson, for $\lambda \in P_\Gamma$ and $w \in W$,
we say that $w$ is $\lambda$-minuscule
if there exists a reduced expression $w=s_{i_1}s_{i_2}\cdots s_{i_t}$
such that
$$s_{i_r}s_{i_{r+1}}\cdots s_{i_t} \lambda = \lambda +
\alpha^{\vee}_{i_r} + \alpha^{\vee}_{i_{r+1}} + \ldots +
\alpha^{\vee}_{i_t}, \, \forall r \in \{1,2,\ldots,t\},$$
where $\alpha_i^\vee \in P_\Gamma$ is the simple coroot corresponding
to $\alpha_i$.  Equivalently (cf.~\cite{stem}), a reduced product
$w = s_1 s_2 \cdots s_{i_t}$ is $\lambda$-minuscule if and only if
$\langle \lambda, \alpha_{i_t}^{\vee} \rangle = -1$ as well as
$\langle s_{i_{r+1}}\cdots s_{i_t} \lambda, \alpha_{i_r}^{\vee}
\rangle = -1$,
for all $r \in \{1,\ldots,t-1\}$, where $\langle \, , \rangle$ is the
Cartan pairing.
\begin{thm}\label{question}
  For $\lambda \in P_{\Gamma}$, there exists a $\lambda$-minuscule
  element $w \in W$ such that $w(\lambda)$ is dominant if and only if
\begin{equation}\label{broer}
  \langle \lambda, \alpha^{\vee} \rangle \geq -1, \, \forall \alpha \in \Delta_+.
\end{equation}
\end{thm}

The proof of this theorem is straightforward, and is given in \S 3.
We also generalize the result to the case of extended Dynkin graphs,
in the following manner.  Let $\widetilde{\Gamma}$ be a simply-laced
extended Dynkin graph, $\widetilde W$ be its Weyl group, and
$R_{\widetilde \Gamma}$ be the root lattice, i.e., the span of the
simple roots $\alpha_i$. Let
$\widetilde{\Delta}_+ \subset R_{\widetilde \Gamma}$ be the set of
positive real roots (i.e., positive-integral combinations $\alpha$ of
simple roots such that $\langle \alpha, \alpha \rangle = 2$).  Define
$P_{\widetilde \Gamma}$ in this case to be the dual to the root
lattice $R_{\widetilde \Gamma}$. Given
$\alpha \in R_{\widetilde \Gamma}$ and
$\lambda \in P_{\widetilde \Gamma}$, denote their pairing by
$\alpha \cdot \lambda$. Let $\delta \in R_{\widetilde \Gamma}$ be the
positive-integral combination of simple roots which generates the
kernel of the Cartan form on $R_{\widetilde \Gamma}$. Finally, for
$\alpha \in \widetilde{\Delta}_+$, let
$\alpha^\vee \in P_{\widetilde \Gamma}$ be the element such that
$\beta \cdot \alpha^\vee = \langle \beta, \alpha \rangle$ for all
$\beta \in \widetilde{\Delta}_+$.  Then, the notion of
$\lambda$-minusculity carries over to this setting.
\begin{thm}\label{question-ext} For nonzero
  $\lambda \in P_{{\widetilde \Gamma}}$, there exists a
  $\lambda$-minuscule element $w \in \widetilde{W}$ such that
  $w(\lambda)$ is dominant if and only if
\begin{enumerate}
\item[(i)]
  $\alpha \cdot \lambda \geq -1, \, \forall \alpha \in
  \widetilde{\Delta}_+$, and
\item[(ii)]
$\delta \cdot \lambda \neq 0$.
\end{enumerate}
\end{thm}
We generalize the theorems above in two directions.  First, we allow
$\lambda$ to be non-integral, i.e., to lie in
$P_{\Gamma} \otimes_{\Z} \R$ (respectively
$P_{\widetilde \Gamma} \otimes_{\Z} \R$) and not just in $P_\Gamma$
(respectively $P_{\widetilde \Gamma}$). Second, we consider all Coxeter
groups, not just finite and affine ones.  For example, in the first
direction, if $\lambda \in P_{\Gamma} \otimes_{\Z} \R$, the notion of
$\lambda$-minuscule Weyl group element should be generalized
accordingly: $w \in W$ is $\lambda$-minuscule if there exists a
reduced expression $w = s_{i_1} \ldots s_{i_t}$ such that
$s_{i_r} \ldots s_{i_t} \lambda = \lambda + \xi_{r}
\alpha_{i_r}^{\vee} + \ldots + \xi_{t} \alpha_{i_t}^{\vee}$
for all $r \in \{1, \ldots, t\}$, for some positive real numbers
$\xi_1, \ldots, \xi_{t} \leq 1$.

In the original situation (for $\lambda \in P_\Gamma$
``integral'' and $\Gamma$ Dynkin), we prove a stronger result:
\begin{thm}\label{question-int}
  Under the assumptions of Theorem \ref{question}, there exists a
  $\lambda$-minuscule element $w \in W$ such that $w(\lambda)$ is
  dominant if and only if
\begin{itemize}
\item[(i)] $\langle \lambda, \alpha_i^{\vee} \rangle \geq -1$ for
  every simple root $\alpha_i$, and
\item[(ii)] For every connected subdiagram $\Gamma' \subseteq \Gamma$,
  the restriction $\lambda|_{\Gamma'}$ is not a negative coroot.
\end{itemize}
\end{thm}
In the theorem, the restriction $\lambda|_{\Gamma'} \in P_{\Gamma'}$
is the unique element such that
$\langle \lambda|_{\Gamma'}, \alpha_i^{\vee} \rangle = \langle \lambda,
\alpha_i^{\vee} \rangle$
for all simple roots $\alpha_i$ associated to the vertices of $\Gamma'$.

We also prove a similar result for extended Dynkin graphs (see Theorem
\ref{intt}), and generalize it so as to include the case where
$\lambda$ lies in a finite Weyl orbit.

\begin{rem}
Condition (\ref{broer}) is equivalent to the
non-negativity of the coefficients of Lusztig's $q$-analogues of
weight multiplicity polynomials (see \cite[Theorem 2.4]{broer}). It is
also equivalent to the vanishing of the higher cohomology groups of
the line bundle that corresponds to $\lambda$ on the cotangent bundle
of the flag variety (op. cit.).  We hope to address and apply this in
future work.
\end{rem}

The paper is organized as follows. The second section introduces the
terminology of Mozes's game of numbers \cite{mozes} and its variant
with a cutoff \cite{qendrim3}, which provides a useful language to
state and prove our results. We also recall some preliminaries on
Dynkin and extended Dynkin graphs. In the third section we solve the
numbers game with a cutoff for Dynkin and extended Dynkin graphs
(Theorem \ref{mt}), in particular proving Theorems \ref{question} and
\ref{question-ext} and the non-integral versions thereof. Next, in \S
4, we give a more explicit solution in the integral case, which proves
Theorem \ref{question-int} and the corresponding result for extended
Dynkin diagrams. In the last section, we generalize Theorem
\ref{question} to the case of arbitrary Coxeter groups.

\subsection{Acknowledgements}
We thank R. Kottwitz for useful comments and M. Boyarchenko for the
opportunity to speak on the topic. The first author is an EPDI fellow
and the second author is an AIM fellow, and both authors were
supported by Clay Liftoff fellowships.  The first author was also
partially supported by the EPSERC Grant EP/F005431/1, and the second
author was partially supported by the University of Chicago's VIGRE
grant. We thank the University of Chicago, MIT, the Max Planck
Institute in Bonn, and the Isaac Newton Institute for Mathematical
Sciences, for hospitality.

\section{The numbers game with and without a cutoff}

In this section we introduce the numbers game with a cutoff, which
provides a useful language to state our results.  We
begin with some preliminaries on Dynkin and extended Dynkin graphs.

\subsection{Preliminaries on Dynkin and extended Dynkin graphs}\label{deds}

We will largely restrict our attention to simply-laced Dynkin and
extended Dynkin graphs.  By this, we mean graphs of type $A_n, D_n,$
or $E_n$, or $\tilde A_n, \tilde D_n$, or $\tilde E_n$. For such a
graph $\Gamma$, let $\Delta$ be the set of (real)\footnote{These are
  sometimes called ``real roots'' in the literature to exclude
  multiples of the so-called imaginary root $\delta$ below, which are
  also roots of the associated Kac-Moody algebra. We will omit the
  adjective ``real.''}
roots of the associated root system, and $\Delta_+$ the set of
positive roots. Let $I$ denote its set of vertices, so that $\alpha_i$
are the simple roots for $i \in I$. Identify $\Z^I$ with the root
lattice (i.e., the integral span of the $\alpha_i$), so that
$\Delta \subseteq \Z^I$, and $\alpha_i \in \Z^I$ are the elementary
vectors. Although we will use subscripts (e.g., $\beta_i$ of
$\beta \in \Z^I$) to denote coordinates, we will never use them for a
vector denoted by $\alpha$, to avoid confusion with the simple roots
$\alpha_i$.

We briefly recall the essential facts about $\Delta_+$ and $\Delta$.
We have $\Delta = \Delta_+ \sqcup (- \Delta_+)$, and
$\Delta_+ = \{ \alpha \in \Z_{\geq 0}^I: \langle \alpha, \alpha
\rangle = 2\}$,
where $\langle \, , \rangle$ is the Cartan form
$$\langle {\alpha_i}, {\alpha_j} \rangle = \begin{cases} 2, & \text{if $i = j$},
  \\ -1, & \text{if $i$ is adjacent to $j$}, \\ 0, &
  \text{otherwise},\end{cases}$$
which is positive-definite in the Dynkin case and
positive-semidefinite in the extended Dynkin case.  It is well known
that $\Delta_+$ is finite in the Dynkin case.  Consider the extended
Dynkin case, and let us switch notation to
$\widetilde{\Gamma}, \widetilde{\Delta}, \widetilde{\Delta}_+$, and
$\widetilde I$.  We may write $\widetilde{\Gamma} \supsetneq \Gamma$
where $\Gamma$ is the Dynkin graph of corresponding type. The vertex
$i_0 = \widetilde{I} \setminus I$ is called an \emph{extending vertex}
(the other extending vertices being obtained as the complements of
different choices of $\Gamma$).  Let $\Delta_+$ the set of positive
roots for $\Gamma$.  There is an inclusion
$\Delta_+ \subset \widetilde{\Delta}_+$ obtained by setting the
coefficient at $i_0$ to zero, and
$\widetilde{\Delta}_+ = (\Delta_+ + \Z_{\geq 0} \delta) \sqcup
(-\Delta_+ + \Z_{>0} \delta)$,
for the unique vector $\delta \in \Z_{>0}^{\widetilde I}$
characterized by $\langle \delta, u \rangle = 0$ for all
$u \in \R^{\widetilde I}$ and $\delta_{i_0} = 1$.

Switching back to $\Gamma, \Delta_+$, and $I$, for either the Dynkin
or extended Dynkin case, we recall the simple reflections.  For any
vertex $i \in I$, let $s_i: \R^I \rightarrow \R^I$ be defined by
$s_i(\beta) = \beta - \langle \beta, {\alpha_i} \rangle
{\alpha_i}$.
It is well known that $\beta \in \Delta_+$ implies
$s_i(\beta) \in \Delta_+$ unless $\beta = {\alpha_i}$, in which case
$s_i({\alpha_i}) = - {\alpha_i}$.  Also, $s_i(\delta) = \delta$ for
all $i$.

For any $\beta \in \Delta_+$, its \emph{height}, $h(\beta)$, is
defined as $h(\beta) = \sum_{i \in I} \beta_i$, where
$\beta = (\beta_i) = \sum_i \beta_i \alpha_i$. Note that $\beta$ may
be obtained from some simple root ${\alpha_i}$ by applying
$h(\beta)-1$ simple reflections, and is not obtainable from any simple
root by applying fewer simple reflections.

\subsection{The numbers game with and without a cutoff}

We first recall Mozes's numbers game \cite{mozes}. Fix an unoriented,
finite graph with no loops and no multiple edges. (For the generalized
version of this game, with multiplicities, see \S \ref{gens}.) Let $I$
be the set of vertices. The \emph{configurations} of the game consist
of vectors $\R^I$. The moves of the game are as follows: For any
vector ${v} \in \R^I$ and any vertex $i \in I$ such that ${v}_i < 0$,
one may perform the following move, called \emph{firing the vertex
  $i$}:
${v}$ is replaced by the new configuration $f_i({v})$, defined by
\begin{equation}
f_i({v})_j =
\begin{cases}
-{v}_i, & \text{if $j = i$}, \\
{v}_j + {v}_i, & \text{if $j$ is adjacent to $i$}, \\
{v}_j, & \text{otherwise}.
\end{cases}
\end{equation}
The entries ${v}_i$ of the vector ${v}$ are called \emph{amplitudes}.
The game terminates if all the amplitudes are nonnegative.  Let us
emphasize that \emph{only negative-amplitude vertices may be
  fired}.\footnote{In
  some of the literature, the opposite convention is used, i.e., only
  positive-amplitude vertices may be fired.}

In \cite{qendrimthesis}, the numbers game \emph{with a cutoff} was
defined: The moves are the same as in the ordinary numbers game, but
the game continues (and in fact starts) only as long as all amplitudes
remain greater than or equal to $-1$.  Such configurations are called
\emph{allowed}. Every configuration which does not have this property
is called \emph{forbidden}, and upon reaching such a configuration the
game terminates (we lose). We call a configuration \emph{winning} if
it is possible, by playing the numbers game with a cutoff, to reach a
configuration with all nonnegative amplitudes.

Call a configuration \emph{losing} if, no matter how the game is
played, one reaches a forbidden configuration.  By definition, any
losing configuration remains so by playing the numbers game. We will
see that the same is true for winning configurations (Theorem
\ref{genmt}).

We now explain how to interpret the results from the introduction in
terms of this language. Let $\Gamma$ be a Dynkin diagram, with set of
vertices $I$. To every element $\lambda \in P_{\Gamma}$ one can
associate naturally an integral configuration of $\Gamma$, still
denoted by $\lambda$, where the amplitude corresponding to the vertex
$\alpha_i$ is given by $\langle \lambda, \alpha_i^{\vee} \rangle$.
Firing the vertex $\alpha_j$ changes these amplitudes to
$\langle s_j(\lambda) , \alpha_i^{\vee} \rangle$, i.e., gives the
natural configuration (on the vertices of $\Gamma$) associated to the
simple reflection $s_j(\lambda)$ of $\lambda$.  In other words, using
the identifications made in the previous subsection between the coroot
space and $\Z^I$, and letting $\cdot$ denote the standard dot product
on $\R^I$, we have
\begin{equation}\label{sifi}
  s_i(\alpha) \cdot v = \alpha \cdot f_i(v), \quad s_i(\alpha) \cdot f_i(v) = \alpha \cdot v,
\end{equation}
for any configuration $v$. In terms of Lie theory, we may think of the
$s_i$ as acting on $\R^I$ with basis given by the simple roots, and
the $f_i$ as acting on the dual $\R^I$, with basis given by the
fundamental coweights. (Formula \eqref{sifi} remains true in the case
of extended Dynkin graphs.)

The existence of an element $w \in W$ such that $w(\lambda)$ is
dominant is then equivalent to the winnability of the usual numbers
game with initial configuration $\lambda$ (and hence, one always
wins). Of course, we want to impose the extra condition that $w$ be
$\lambda$-minuscule, which is equivalent to imposing the $-1$ cutoff
to the numbers game. Thus, Theorem \ref{question} gives a
characterization of the winning configurations ${v} \in \Z^I$ for the
numbers game with a cutoff, where
$v_i = \langle \lambda, \alpha_i^{\vee} \rangle$,
$\lambda \in P_{\Gamma}$, and the graph $\Gamma$ is a Dynkin
one. Later on, we will give similar descriptions in terms of the
numbers game with a cutoff for the other results stated in the
introduction. 

Note that in the paragraph above we only considered the case of
integral $\lambda$, but the analogy holds in the non-integral case as
well, and now we study the winnability of the numbers game with a
cutoff with real amplitudes, where we may fire any vertex with
amplitudes from $[-1,0)$ and not just those with amplitude $-1$ as in
the integral case.

The language of the numbers game with a cutoff is useful because it
makes apparent certain phenomena that already occur without the bound
of 1 or indeed with a different bound. It also allows one to use
results from the usual Mozes's numbers game, which has been widely
studied (cf. \cite{Pro-bru, Pro-min, DE, Erik-no1, Erik-no2, erikconf,
  Erik-no3, Erik-no4, eriksson, Wild-no1, Wild-no2}),\footnote{Mozes's
  numbers game originated from (and generalizes) a 1986 IMO problem.}
and yields useful algorithms for computing with the root systems and
reflection representations of Coxeter groups (see \cite[\S 4.3]{BB}
for a brief summary).

Finally, we recall some basic results about the usual numbers game,
and why it exhibits special behavior in the Dynkin and extended Dynkin
cases:
\begin{prop}
\begin{enumerate}
\item[(i)] \cite{mozes} If the usual numbers game terminates, then it
  must terminate in the same number of moves and at the same
  configuration regardless of how it is played.
\item[(ii)] In the Dynkin case, the usual numbers game must terminate.
\item[(iii)] \cite{erikconf} In the extended Dynkin case, the usual numbers game terminates if and only if $\delta \cdot v > 0$.
\item[(iv)] \cite{erikconf} Whenever the usual numbers game does not
  terminate, it reaches infinitely many distinct configurations,
  except for the case of an extended Dynkin graph where
  $\delta \cdot v = 0$, in which case only finitely many
  configurations are reached \emph{(}i.e., the game
  ``loops''\emph{)}.\footnote{Stronger results were stated in
    \cite{erikconf}, and a detailed study appears in \cite{GSS}.}
\end{enumerate}
\end{prop}
Thus, provided we can determine which configurations are winning (for
the numbers game with a cutoff) in the Dynkin case and the extended
Dynkin case, then with the additional condition $\delta \cdot v> 0$,
these are also the ones that terminate in a nonnegative configuration,
and this configuration (and the number of moves required to get there)
is unique.

\section{The (extended) Dynkin case}
\begin{thm} \label{mt}
In the Dynkin case, a configuration ${v}$ is winning if and only if
\begin{equation} \label{rtconddyn}
  \alpha \cdot {v} \geq -1, \quad \forall \alpha \in \Delta_+.
\end{equation}
Otherwise, ${v}$ is losing.

In the extended Dynkin case, ${v} \neq 0$ is winning if and only if
both
\begin{equation} \label{rtconded}
  \alpha \cdot {v} \geq -1, \quad \forall \alpha \in \widetilde{\Delta}_+,
\end{equation}
and $\delta \cdot {v} \neq 0$.  If \eqref{rtconded} is satisfied but
$\delta \cdot {v} = 0$ \emph{(}and ${v} \neq 0$\emph{)}, then ${v}$ is
looping and the game cannot terminate.  Finally, if \eqref{rtconded}
is not satisfied \emph{(}e.g., if $\delta \cdot {v} < 0$\emph{)}, then
${v}$ is losing.
\end{thm}

\begin{rem}
  Theorem \ref{mt} implies Theorems \ref{question} and
  \ref{question-ext}, as well as their ``non-integral'' versions.
\end{rem}

The above theorem shows, in particular, that exactly one of the
following is true: ${v}$ is winning, looping, or losing.

To prove the theorem, it is helpful to introduce the set
\begin{equation}
  X_{{v}} := \{(\alpha, \alpha \cdot {v}) \mid \alpha \in \Delta_+, \alpha \cdot {v} < 0 \}.
\end{equation}

Consider the projections
\begin{equation}
\xymatrix{
& X_{{v}} \ar@{^{ (}->}[ld]_{\pi_1} \ar[rd]^{\pi_2} & \\ \Delta_+ & & \R_{<0}.
}
\end{equation}
Each time a vertex, say $i \in I$, is fired, there is a natural
isomorphism $X_{{v}} \setminus \{({\alpha_i}, v_i)\} \iso X_{{f_i
    v}}$,
with
$(\alpha, \alpha \cdot {v}) \mapsto (s_i \alpha, \alpha \cdot {v}) =
(s_i \alpha, s_i \alpha \cdot f_i {v})$.
The set $X_v$ is defined similarly in the extended Dynkin case, with
$\Delta_+$ replaced by $\widetilde{\Delta}_+$, and there is still a
natural isomorphism
$X_{{v}} \setminus \{({\alpha_i}, v_i)\} \iso X_{{f_i v}}$.

\begin{proof} In the Dynkin case, $X_{{v}}$ is finite. Since the size
  decreases by one in each step, removing an element whose second
  projection is the amplitude at the vertex which is fired, we see
  that the game is won precisely when $\pi_2(X_{{v}}) \subset [-1,0)$,
  and otherwise it is lost. The former is equivalent to
  \eqref{rtconddyn}.

  In the extended Dynkin case, the game is won precisely when
  $X_{{v}}$ is finite and $\pi_2(X_{{v}}) \subset [-1,0)$; finiteness
  is equivalent to $\delta \cdot {v} > 0$.  The condition
  $\pi_2(X_{{v}}) \subset [-1,0)$ is equivalent to \eqref{rtconded},
  and implies $\delta \cdot {v} \geq 0$, so for ${v}$ to be winning we
  only need to additionally assume that $\delta \cdot {v} \neq 0$.

Since, in the extended Dynkin case, a game that is not won is either
lost or loops, it remains to show that ${v}$ is losing precisely when
there exists $\alpha \in \widetilde{\Delta}_+$ with
$\alpha \cdot {v} < -1$, i.e., when
$\pi_2(X_{{v}}) \not \subset [-1,0)$.  It is clear that the condition
is required for ${v}$ to be losing.  Thus, suppose that
$\alpha \cdot {v} < -1$ for some $\alpha \in \widetilde{\Delta}_+$. We
will show that ${v}$ is losing.  We induct on the height of $\alpha$.
Suppose $v_i < 0$, and that we fire the vertex $i$. Consider two
cases: first, suppose that $h(s_i \alpha) < h(\alpha)$.  Then,
$s_i \alpha \cdot f_i {v} < -1$ and $h(s_i \alpha) < h(\alpha)$,
completing the induction.  Next, suppose
$h(s_i \alpha) \geq h(\alpha)$, i.e., $s_i \alpha - \alpha$ is a
nonnegative multiple of ${\alpha_i}$.  Then,
$\alpha \cdot f_i {v} \leq s_i \alpha \cdot f_i {v}$ (since
$(f_i {v})_i > 0$), and $s_i \alpha \cdot f_i {v} = \alpha \cdot {v}$.
Thus, we may leave $\alpha$ unchanged.  If we eventually fire a vertex
$i \in \widetilde{I}$ such that $h(s_i \alpha) < h(\alpha)$, the
induction is complete. Otherwise, we would be playing the game only on
a Dynkin subgraph, which would have to terminate in finitely many
moves, and therefore reach a forbidden configuration (since
$\pi_2(X_{{v}}) \not \subset [-1,0)$).
\end{proof}

Note that only finitely many inequalities in \eqref{rtconded} are
required: since \eqref{rtconded} implies $\delta \cdot v \geq 0$,
\eqref{rtconded} is equivalent to the conditions
$\delta \cdot {v} \geq 0$, $\alpha \cdot v \geq -1$, and
$(\delta - \alpha) \cdot v \geq -1$ for all $\alpha$ which are
positive roots of a corresponding Dynkin subgraph obtained by removing
an extending vertex. So, it is enough to assume \eqref{rtconded} for
$\alpha \in \Delta_+ \cup (\delta - \Delta_+)$, which is finite.

\begin{cor} \label{dzcor} If $\delta \cdot {v} = 0$, then the game
  loops \emph{(}and cannot terminate\emph{)} if and only if, after
  removing an extending vertex, both ${v}$ and $-{v}$ are winning.
\end{cor}
\begin{proof}
  This follows from the fact that
  $\widetilde{\Delta}_+ = (\Delta_+ + \Z_{\geq 0} \delta) \sqcup
  (-\Delta_+ + \Z_{> 0} \delta)$.
\end{proof}
Another interpretation of the above corollary is the following: ${v}$
continues indefinitely if and only if the restriction of ${v}$ to the
complement of an extending vertex cannot reach a forbidden
configuration by playing numbers game forwards \emph{or backwards}
(i.e., firing vertices with positive instead of negative amplitudes).

\begin{rem} T. Haines pointed out that Theorem \ref{mt} implies
  \cite[Lemma 3.1]{Haines}: for every dominant
  minuscule\footnote{Recall that minuscule means that
    $\langle \mu, \alpha \rangle \in \{-1, 0, 1\}$ for all
    $\alpha \in \Delta$.}
  coweight $\mu$ and every coweight $\lambda \in W \mu$, there exists
  a sequence of simple roots $\alpha_1, \ldots, \alpha_p$, such that
  $s_1(\mu) = \mu-\alpha_1^\vee$,
  $s_2 s_1 \mu = \mu - \alpha_1^\vee - \alpha_2^\vee, \ldots,$ and
  $\lambda = s_p s_{p-1} \cdots s_1(\mu) = \mu - \alpha_1^\vee -
  \cdots - \alpha_p^\vee$.
\end{rem}

\section{The integral case}\label{ints}

Of particular relevance is the case of integral configurations
${v} \in \Z^I$.  Below, we apply Theorem \ref{mt} to give a
surprisingly simple, explicit description of the losing and looping
integral configurations in the Dynkin and extended Dynkin cases.

To state the theorem, we will make use of the interpretation of
configurations ${v} \in \R^I$ as coweights. In particular, as in the
introduction, for every Dynkin graph $\Gamma$, and every root
$\alpha \in \Delta_+$, there is an associated coroot configuration
$\alpha^{\vee} \in \Z^I$, in the basis of fundamental coweights,
uniquely defined by
$\beta \cdot \alpha^{\vee} = \langle \beta, \alpha \rangle$ for all
$\beta$, using the Cartan form as in \S \ref{deds}. For every extended
Dynkin graph $\widetilde{\Gamma}$, Dynkin subgraph $\Gamma$, and
$\alpha \in \widetilde{\Delta}_+$, we also have the configuration
$\alpha^{\vee}$ defined in the same way; in particular,
$\delta \cdot \alpha^{\vee} = 0$ (and the $\alpha_i^\vee$ are linearly
dependent). Let $\omega_i \in \Z^I$ be the elementary vector, viewed
as a configuration (i.e., in the Dynkin case, the $i$-th fundamental
coweight).\footnote{We use distinct notation $\alpha_i, \omega_i$ for
  the same vector in $\Z^I$ depending on whether it is viewed as a
  simple root or a configuration, to avoid confusion.}
Thus, $\alpha_i \cdot \omega_j = \delta_{ij}$.  For
$\beta \in \Delta_+$ or $\widetilde{\Delta}_+$, let its
\emph{support}, $\text{supp}(\beta)$, be the (connected) subgraph on
which its coordinates $\beta_i$ are nonzero.
\begin{thm} \label{intt}
\begin{itemize}
\item[(i)] An integral configuration $v$ on a Dynkin graph is winning
  if and only if
\begin{itemize}
\item[(1)] $v_i \geq -1$ for all $i$, and
\item[(2)] For all $\alpha \in \Delta_+$,
  $v|_{supp(\alpha)} \neq -\alpha^{\vee}$;
\end{itemize}
\item[(ii)] An integral configuration $v$ on an extended Dynkin graph
  is winning if and only if \emph{(1)} and \emph{(2)} are satisfied
  \emph{(}with $\alpha \in \widetilde \Delta_+$\emph{)}, and
  furthermore,
\begin{itemize}
\item[(3)] $v \neq -\omega_i$ for any extending vertex $i$.
\end{itemize}
\item[(iii)] An integral configuration on an extended Dynkin graph is
  looping if and only if it is in the Weyl orbit of a vector
  $\mu= \omega_i - \omega_{i'}$ for distinct extending vertices
  $i, i'$.  In this case, the numbers game can take the configuration
  to and from such a vector $\mu$.
\end{itemize}
\end{thm}

\begin{rem}
  The above result implies Theorem \ref{question-int}, as well as the
  extended Dynkin version thereof.
\end{rem}
As in the introduction, for $\Gamma' \subseteq \Gamma$, with vertex
sets $I' \subseteq I$, the restriction $v|_{\Gamma'}$ is the
restriction $\R^I \onto \R^{I'}$ of coordinates.

We remark that an alternative way to state parts (i) and (ii) above is
that the losing configurations on (extended) Dynkin diagrams which are
winning on all proper subgraphs, which we call the \emph{minimal
  losing configurations},
are exactly those of the form $-\beta^{\vee}$ for fully supported
roots $\beta$, which in the extended Dynkin case also satisfy
$\beta_i \leq \delta_i$ for all $i$, and $-\omega_j$ for extending
vertices $j$, together with the one-vertex forbidden configurations.

Here, we have used that $(\beta + c \delta)^\vee = \beta^\vee$ for
all $c \in \Z$, so that in part (ii) it suffices to assume that
$\beta \in \widetilde{\Delta}_+$ satisfies $\beta_i \leq \delta_i$
for all $i$, i.e., $\beta_i \leq 1$ for all extending vertices $i$.
In fact, we can further restrict to the case of roots $\beta$ that
are supported on a Dynkin subdiagram, in exchange for adding the
condition that $v_{supp(\gamma)} \neq \gamma^{\vee}$ for all positive
roots $\gamma$ such that $\gamma_i = 0$ at all extending vertices
$i$. This is because the fully supported roots $\beta$ such that
$\beta_i \leq \delta_i$ for all $i$ are exactly $\delta - \gamma$
where $\gamma \in \widetilde{\Delta}_+$ satisfies $\gamma_i = 0$ at
all extending vertices, and then $-\beta^\vee = \gamma^\vee$.

As a special case of (ii), for $\widetilde{A_n}$ (with $n \geq 1$),
the only integral losing configurations which are winning on all
proper subgraphs are $- {\omega_i}$ for all $i$.  Also, by (iii),
there is no looping integral configuration on $\widetilde{E_8}$ (but
these exist for all other extended Dynkin graphs).

\begin{proof}
  (i) Following the discussion above, we show that the minimal losing
  configurations on Dynkin graphs with more than one vertex are
  exactly $-\beta^\vee$ for fully supported $\beta \in \Delta_+$.
  Note that it is clear that such configurations are minimal losing
  configurations, since $\beta \cdot (-\beta^\vee) = -2$ and
  $\gamma \cdot (-\beta^\vee) \in \{-1,0,1\}$ for all
  $\gamma \in \Delta_+ \setminus \{\beta\}$. Thus, we only need to
  show that there are no other minimal losing configurations (other
  than one-vertex ones).

  For any minimal losing configuration $v \in \Z^I$, Theorem \ref{mt}
  implies the existence of $\beta \in \Delta_+$ such that
  $\beta \cdot v \leq -2$.  By minimality, all such $\beta$ are
  fully supported.  It suffices to prove that, when $\beta$ is not
  simple (i.e., the graph has more than one vertex),
  $v = -\beta^\vee$.  We prove this by induction on the height of
  $\beta$, considering all Dynkin graphs simultaneously.

  Let $i$ be a vertex such that $h(s_i \beta) < \beta$, i.e.,
  $\langle \beta, {\alpha_i} \rangle = 1$.  It follows that
  ${v}_i = -1$; otherwise, $s_i \beta \cdot {v} \leq -2$, a
  contradiction. 
  Since $s_i \beta \cdot f_i v \leq -2$, we deduce from the inductive
  hypothesis that the restriction of $f_i v$ to the support of
  $s_i \beta$ coincides with $-(s_i \beta)^\vee$. Since
  $-((s_i \beta)^{\vee})_i = (\beta^\vee)_i = 1$, we deduce that
  $f_i v = -(s_i \beta)^\vee$ and hence $v = -\beta^\vee$, as desired.

  (ii) We prove that the minimal losing configurations in the extended
  Dynkin case are exactly $-\beta^\vee$ for fully supported
  $\beta \in \widetilde \Delta_+$ satisfying $\beta_i \leq \delta_i$
  for all $i$, and $-\omega_i$ for extending
  vertices $i$.  The former configuration is a minimal losing
  configuration by the same argument as in the Dynkin case, and
  $-\omega_i$ is a minimal losing configuration since
  $\delta \cdot -\omega_i = -1 < 0$ (so $-\omega_i$ is losing) and
  $\beta \cdot -\omega_i = -\beta_i \in \{-1, 0\}$ for all
  $\beta \in \widetilde \Delta_+$ (so $-\omega_i$ is winning on all
  Dynkin subdiagrams).  Hence, it suffices to prove that there are no
  other minimal losing configurations.

  Let ${v}$ be an integral losing configuration which is winning on
  all proper subdiagrams, and let $\beta \in \widetilde{\Delta}_+$ be
  of minimal height such that $\beta \cdot {v} \leq -2$.  Once again,
  we can induct on the height of $\beta$.  We reach the desired
  conclusion unless $\beta = c \delta + {\alpha_i}$ for some
  $c \geq 1$ and $i \in \widetilde I$, so assume this.
  Since $v_i \geq -1$, it follows that $\delta \cdot {v} \leq -1$.
  Moreover, fix an associated Dynkin subdiagram $\Gamma$. Then, for
  all $\gamma \in \Delta_+$, we must have
  $\gamma \cdot {v} \in \{-1,0\}$ (since
  $(\delta - \gamma) \cdot {v} \geq -1$ and $\gamma \cdot {v} \geq -1$
  by minimality of $\beta$).  In particular, ${v}_j \in \{-1,0\}$ for
  all $j$.  In this case, in order not to be losing on a Dynkin
  subdiagram, we must have ${v} = - {\omega_i}$, where $i$ is an
  extending vertex.

  (iii) Let $i$ be an extending vertex, and let ${v} \in \Z^I$ satisfy
  $\delta \cdot {v} = 0$ but ${v} \neq 0$. If we play the numbers game
  by firing only vertices other than $i$, we must eventually obtain
  either a forbidden configuration (if the restriction of ${v}$ to the
  complement of $i$ is losing) or a configuration whose sole negative
  amplitude occurs at $i$.  In the latter case, in order to not be
  forbidden, we must have $-1$ at the vertex $i$, and hence, in order
  to satisfy $\delta \cdot {v}= 0$, there can only be one positive
  amplitude, it must be $1$, and it must occur at another extending
  vertex, say $i'$.  So, ${v}$ is winning when restricted to the
  complement of $i$ if and only if one can obtain
  $\mu = {\omega_{i'}} - \omega_{i}$ from ${v}$.  This implies that
  $v$ is in the same Weyl orbit as $\mu$.  On the other hand, if $v$
  is in the Weyl orbit of $\mu$, then $\delta \cdot v = 0$ and the
  usual numbers game loops, and since
  $\alpha \cdot v \in \{-1, 0, 1\}$ for all
  $\alpha \in \widetilde{\Delta}_+$, the numbers game with a cutoff
  also loops.  Hence, the conditions that $v$ is looping, that $v$ is
  in the Weyl orbit of such a $\mu$, and that $\mu$ can be obtained
  from $v$ by playing the numbers game with a cutoff, are all
  equivalent.  Since, in this case, $-v$ is also looping, we see also
  that $-v$ can reach a configuration $\nu = \omega_j - \omega_{i'}$
  for some extending vertex $j$, and since $\nu$ is in the same Weyl
  orbit as $-\mu$, we must have $\nu = - \mu$ (since $-\mu$ and $\nu$
  are dominant on the complement of $i'$).  Hence, $v$ can be obtained
  from $\mu$ by playing the numbers game, which proves the remainder
  of the final assertion.
\end{proof}

\begin{rem} In the Dynkin case, the above may be interpreted as saying
  that every losing integral configuration which is winning on all
  proper subgraphs is obtainable from the maximally negative coroot by
  playing the numbers game: this configuration is the one with
  ${v}_i = -1$ when $i$ is adjacent to the extending vertex of
  $\widetilde{\Gamma}$, and ${v}_i = 0$ otherwise.  On the other hand,
  in the non-integral case, losing configurations are not necessarily
  obtainable from nonpositive ones by playing the numbers game: for
  example, on $D_4$, one may place $-1$ at all three
  endpoint vertices, and $\frac{3}{2}$ at the node.
\end{rem}

\begin{rem} Note that the extended Dynkin case with
  $\delta \cdot {v} \geq 0$ and ${v}$ losing, integral, and winning on
  all subgraphs may similarly be described as those configurations
  obtainable from
  ${\alpha_i}^\vee = 2 {\omega_i} - \sum_{j \text{ adjacent to } i}
  {\omega_j}$,
  for $i$ not an extending vertex, by playing the numbers game.  This
  contrasts with the nonintegral case: see the next remark.
\end{rem}

\begin{rem}
  In the extended Dynkin case, it is perhaps surprising that all
  losing integral configurations with $\delta \cdot {v} > 0$ are also
  losing on a proper subgraph. This is not true in the non-integral
  case (except in the case $\widetilde{A_n}$): e.g., one may take a
  configuration $\beta^{\vee} + \varepsilon {\omega_i}$, for
  $\beta \in \widetilde{\Delta}_+$ which satisfies $\beta_j = 0$ for
  all extending vertices $j$, and
  $\varepsilon \in (0,\frac{1}{\delta_i})$ for any fixed
  $i \in \widetilde{I}$. Similarly, one may find losing configurations
  with $\delta \cdot {v} = 0$ which are winning on all Dynkin
  subgraphs, but are not $\beta^{\vee}$ for $\beta \in \Delta_+$
  (although there are still none for $\widetilde{A_n}$): for example,
  $\varepsilon \beta^{\vee}$ for $\varepsilon \in (\frac{1}{2}, 1)$
  and $\beta$ as before.  For another example, we can take any
  configuration in $\widetilde{D_n}$ with values $a, b, c, d \geq -1$
  at exterior vertices such that $\sigma := \frac{a+b+c+d}{2} < -1$
  and $\sigma - x \geq -1$ for all $x \in \{a,b,c,d\}$.  Finally,
  there are many more losing nonintegral configurations with
  $\delta \cdot {v} < 0$ that are winning on all subgraphs than just
  $-{\omega_i}$ for $i$ an extending vertex: for example,
  $-{\omega_i} + u$ for any nonnegative vector $u$ such that
  $\delta \cdot u < 1$.
\end{rem}

\section{Generalization to arbitrary graphs with multiplicities} \label{gens}

In \cite{mozes, eriksson}, the numbers game was stated in greater
generality than the above.  Namely, in addition to a graph with vertex
set $I$ (and no loops or multiple edges), we are given a Coxeter group
$W$ with generators $s_i, i \in I$ and relations $(s_i s_j)^{n_{ij}}$
for $n_{ij} \in \{1, 2, \ldots\} \cup \{\infty\}$, together with a
Cartan matrix $C = (c_{ij})_{i,j \in I}$, such that $c_{ii} = 2$ for
all $i$, $c_{ij} = 0$ whenever $i$ and $j$ are not adjacent, and
otherwise $c_{ij}, c_{ji} < 0$ and either
$c_{ij} c_{ji} = 4 \cos^2(\frac{\pi}{n_{ij}})$ (when $n_{ij}$ is
finite) or $c_{ij} c_{ji} \geq 4$ (when $n_{ij} = \infty$).

We recall that the numbers game is modified as follows in terms of
$C$: The configurations are again of the form ${v} \in \R^I$, and, we
may fire the vertex $i$ in a configuration ${v} \in \R^I$ if and only
if the amplitude ${v}_i < 0$. The difference is that the new
configuration $f_i({v}_i)$ is now given by
\begin{equation}
f_i({v})_j = v_j - c_{ij} v_i.
\end{equation}
We call this the \emph{weighted} numbers game.  The non-weighted numbers game is recovered in the case $c_{ij} =-1$ for all adjacent $i,j$.

The standard reflection action of $W$ on $\R^I$ is given by
\begin{equation}
s_i(\beta)_j =
\begin{cases} \beta_j, & \text{if $j \neq i$}, \\
-\beta_i - \sum_{k \neq i} c_{ik} \beta_k, & \text{if $j = i$}.
\end{cases}
\end{equation}
Recall from \cite{eriksson} that, in this situation, the usual numbers
game is \emph{strongly convergent}: if the game can terminate, then it
must terminate, and in exactly the same number of moves and arriving
at the same configuration, regardless of the choices made.

We remark that, while it is standard to take $C$ to be symmetric,
there are cases when this is not desired, particularly for the
non-simply-laced Dynkin diagrams $\Gamma$, where $C$ can be taken to
be integral only if allowed to be non-symmetric.  In these cases, if
we choose $C$ to be integral, playing the numbers game on $\Gamma$ is
equivalent to playing the numbers game without multiplicities on a
simply-laced diagram $\Gamma'$ with some symmetry group $S$, such that
$\Gamma' / S = \Gamma$, if we restrict to $S$-invariant configurations
on $\Gamma'$, where we allow simultaneous firing of any orbit of
vertices under $S$ (since these orbits consist of nonadjacent
vertices, it makes sense to fire them simultaneously).

Let $\Delta = \bigcup_{i \in I} W {\alpha_i}$ be the set of
(\emph{real}) \emph{roots}.\footnote{Note that, when the Cartan matrix
  $C$ is associated to a nonreduced root system (i.e., $BC_n$), then
  $\Delta$ is a proper subset of the whole root system, which does not
  contain $2\alpha$, for any simple root $\alpha$.}
Let $\Delta_+ \subset \Delta$ be the subset of \emph{positive roots}:
these are the elements whose entries are nonnegative. Note that, by a
standard result (see \cite[Proposition 4.2.5]{BB}),
$\Delta = \Delta_+ \sqcup (- \Delta_+)$.

Finally, we recall a useful partial ordering from, e.g., \cite[\S
4.6]{BB}.  For $\beta \in \Delta_+$, we say that $\beta < s_i \beta$
if and only if $\beta_i < (s_i \beta)_i$. Generally, for
$\alpha, \beta < \Delta_+$, we say $\alpha < \beta$ if there exists a
sequence
$\alpha < s_{i_1} \alpha < s_{i_2} s_{i_1} \alpha < \cdots < s_{i_m}
s_{i_{m-1}} \cdots s_{i_1} \alpha = \beta$.
The argument of \cite[Lemma 4.6.2]{BB} shows that this is a graded
partial ordering. The grading, $\dep(\alpha)$, called the
\emph{depth}, is defined to be the minimum number of simple
reflections required to take $\alpha$ to a negative root.  Thus,
$\alpha < s_i \alpha$ implies $\dep(s_i \alpha) = \dep(\alpha) + 1$.

\begin{thm}\label{genmt} Let $\Gamma, C$ be associated to a Coxeter
  group. Assume that $C$ satisfies $c_{ij} = c_{ji}$ whenever $n_{ij}$
  is odd \emph{(}and finite\emph{)}. Then, ${v}$ can reach a forbidden
  configuration if and only if $\beta \cdot {v} < -1$ for some
  $\beta \in \Delta_+$, and in this case, the minimum number of moves
  required to take ${v}$ to a forbidden configuration is
\begin{equation}
  m({v}) := \text{min} \{\dep(\beta)-1 \mid \beta \cdot {v} < -1, \beta \in \Delta_+\}.
\end{equation}
Furthermore, if ${v}_i < 0$, then $m(f_i {v}) \in \{m({v}), m({v})-1\}$.
\end{thm}
Note that, in the non-simply-laced Dynkin cases with $C$ integral, we
may always take $c_{ij} = c_{ji}$ whenever $n_{ij}$ is odd (and in
these cases, this implies $n_{ij} = 3$), so the theorem applies.  
\begin{cor} Under the assumptions of the theorem, ${v}$ is winning if and only if the usual numbers game terminates and
\begin{equation} \label{rtcond2}
\alpha \cdot {v} \geq -1, \forall \alpha \in \Delta_+.
\end{equation}
Moreover, if \eqref{rtcond2} is not satisfied and the usual numbers
game terminates, then ${v}$ is losing.
\end{cor}
Also, under the hypotheses of the theorem, any winning configuration
remains so regardless of what moves are made.

We can also make a statement for arbitrary $C$ and $\Gamma$:

\begin{thm}\label{onedirthm} If $C$ and $\Gamma$ are arbitrary
  \emph{(}associated to a Coxeter group\emph{)}, then ${v}$ can reach
  a forbidden configuration if and only if there exists
  $\beta \in \Delta_+$ and $i \in I$ such that both
  $\beta \cdot {v} < -1$ and $\beta > {\alpha_i}$.
In this case, the minimum number of moves required to reach a forbidden configuration is
\begin{equation} \label{mpdfn} m'({v}) := \text{min} \{\dep(\beta)-1
  \mid \beta \cdot {v} < -1, \text{ and there exists $i \in I$ with }
  \beta > {\alpha_i}\}.
\end{equation}
Moreover, in this case, if $i \in I$ is such that ${v}_i < 0$, then
$m'(f_i {v}) \geq m'({v}) - 1$ \emph{(}provided $m'(f_i {v})$ is
defined, i.e., $f_i {v}$ can reach a forbidden configuration\emph{)}.
\end{thm}

The difference from Theorem \ref{genmt} is that we added the condition
$\beta > {\alpha_i}$, and replaced the equality for $m$ under numbers
game moves by an inequality.

We remark that the usual numbers game, beginning with ${v}$,
terminates if and only if
\begin{equation}
\# \mathbb{P}\{\beta \in \Delta_+ \mid \beta \cdot {v} < 0\} < \infty,
\end{equation}
for arbitrary $\Gamma, C$, where $\mathbb{P}$ means modding by scalar
multiples, since each move decreases the size of this set by one. (We
do not need to mod by scalar multiples if $c_{ij} = c_{ji}$ whenever
$n_{ij}$ is odd.)  So, this gives a completely root-theoretic
description of the winning conditions above.\footnote{Also, this
  observation easily implies the main results (Theorems 2.1 and 4.1)
  of \cite{DE}: if ${v}_i \leq 0$ for all $i$ and ${v} \neq 0$, then
  the usual numbers game can only terminate if $\Gamma, C$ are
  associated to a finite Coxeter group: otherwise (assuming $\Gamma$
  is connected), infinitely many elements $\beta \in \Delta_+$ which
  are not multiples of each other satisfy $\beta \cdot {v} < 0$: note
  that, for each $i \in I$, the set $\mathbb{P}(W {\alpha_i})$
  essentially does not depend on the choice of $C$ for a given Coxeter
  group.}

For the finite and affine cases, we have the following corollary,
which generalizes Theorem \ref{mt}.  As before, in the affine case,
let $\delta \in \R_{> 0}^I$ be the additive generator of the semigroup
$\{\delta' \in \R_{> 0}^I \mid \alpha \in \Delta_+ \Rightarrow
\alpha+\delta' \in \Delta_+\}$.
In particular, $\langle \delta, \alpha \rangle = 0$ for all
$\alpha \in \Delta$.
\begin{cor} \label{genmtcor} Let $\Gamma, C$ be associated to a finite
  or affine Coxeter group and let ${v}$ be a nonzero configuration.
  Then, exactly one of the following is true:
\begin{enumerate}
\item[(a)] \eqref{rtcond2} is satisfied, and $\delta \cdot {v} \neq
  0$:
  then ${v}$ is winning, and cannot reach a forbidden configuration.
\item[(b)] \eqref{rtcond2} is satisfied but $\delta \cdot {v} = 0$:
  then ${v}$ is looping, and cannot reach a forbidden configuration.
\item[(c)] \eqref{rtcond2} is not satisfied.  Then, provided
  $c_{ij} = c_{ji}$ whenever $n_{ij}$ is odd, ${v}$ is losing.
\end{enumerate}
\end{cor}

Note that, by Theorem \ref{onedirthm}, we can strengthen this slightly
by replacing \eqref{rtcond2} by the condition that
$\alpha \cdot {v} \geq -1$ only for $\alpha$ such that
$\alpha > {\alpha_i}$ for some $i \in I$.

\begin{proof}[Proof of Corollary \ref{genmtcor}] (a) In the affine
  case, $\delta \cdot {v} > 0$, so in either case, the usual numbers
  game terminates. Then, ${v}$ is winning by Theorem \ref{onedirthm},
  and a forbidden configuration cannot be reached.

  (b) ${v}$ is looping, as in the simply-laced case, since the usual
  numbers game cannot terminate, and the configuration is uniquely
  determined by its restriction to a subgraph obtained by removing an
  extending vertex, where the configuration remains in the orbit of
  the restriction of ${v}$ under the associated finite Coxeter
  group. The rest follows from Theorem \ref{onedirthm}.

(c) In this case (we assume $c_{ij} = c_{ji}$ whenever $n_{ij}$ is
odd), ${v}$ can reach a forbidden configuration. Moreover, in the
proof of Theorem \ref{genmt}, we see that there always exists a vertex
$i \in I$ so that, for any configuration ${v}'$ obtained from ${v}$ by
firing vertices other than $i$, we have $m(f_i {v}') = m({v}') - 1$.
In the affine Coxeter group case, in order for the numbers game to
continue indefinitely, all vertices must be fired infinitely many
times. This proves the result.
\end{proof}

\begin{rem} 
  The weakened conclusions of Theorem \ref{onedirthm} are
  needed. Indeed, if $c_{ij} \neq c_{ji}$ for some $i,j$ with $n_{ij}$
  odd, then it is possible that a winning configuration can become a
  losing one. For example, take $I = \{1,2\}$ and
  $C = \begin{pmatrix} 2 & -2 \\ -\frac{1}{2} & 2 \end{pmatrix}$, with
  $n_{12} = 3$.  Then, the configuration
  $(-\frac{1}{2}, -\frac{1}{2})$ is winning under the sequence
  $(-\frac{1}{2}, -\frac{1}{2}) \mapsto (-\frac{3}{4}, \frac{1}{2})
  \mapsto (\frac{3}{4}, -1) \mapsto (\frac{1}{2}, 1)$,
  but if we instead fired vertex $1$ first, we would get
  $(\frac{1}{2}, -\frac{3}{2})$, which is forbidden.
\end{rem}

\begin{rem} It is natural to ask what can happen in the numbers game
  with a cutoff if it continues indefinitely.  Suppose this happens
  and that $\Gamma'$ is the subgraph on vertices which are fired
  infinitely many times.  If $\Gamma'$ corresponds to an affine
  Coxeter group, then the configuration restricted to $\Gamma'$ is
  looping, and in this case, in order for a forbidden configuration
  not to be reached, $\Gamma'$ must be the whole graph (assuming that
  our whole graph is connected).  Otherwise, if our graph is not
  affine, then $\Gamma'$ cannot be associated to an affine or finite
  Coxeter group.  Then, for any affine subgraph
  $\Gamma_0 \subseteq \Gamma'$ (where by this we allow reducing the
  numbers $n_{ij}$ for edges between vertices of $\Gamma_0$), the
  inner product of the restriction of ${v}$ with the associated
  $\delta_0$ must remain positive, and the value must be decreasing.
  It must converge to some nonnegative number, and hence all
  amplitudes of vertices in $\Gamma'$ must converge to zero.  In
  particular, the configuration $v$ must converge to some limiting
  allowed configuration (which is zero on $\Gamma'$), and one could
  continue the numbers game from this limit if desired.  Note that, in
  the case that $c_{ij} = c_{ji}$ for all odd $n_{ij}$, we must also
  have $\alpha \cdot {v} > -1$ for all $\alpha \in \Delta_+$ supported
  on $\Gamma'$, i.e., ${v}|_{\Gamma'}$ cannot reach a forbidden
  configuration by playing the numbers game on $\Gamma'$.
\end{rem}

\subsection{Proof of Theorems \ref{genmt} and \ref{onedirthm}}
We will use the following lemma which is interesting in itself (and is
the connection between the two theorems):
\begin{lemma} \label{curlem} If $\Gamma, C$ are such that
  $c_{ij} = c_{ji}$ whenever $n_{ij}$ is odd, then for all
  $\beta \in \Delta_+$, we have ${\alpha_i} \leq \beta$ for some
  $i \in I$.
\end{lemma}

We remark that it is well known (and obvious) that the lemma holds
when $C$ is symmetric.

\begin{proof}
  The case $n_{ij}$ is odd is exactly the case when, on the subgraph
  with vertices $i$ and $j$ only, ${\alpha_i}$ is in the $W$-orbit of
  some positive multiple of ${\alpha_j}$ and vice-versa (and this
  multiple is $1$ if and only if $c_{ij} = c_{ji}$).  Thus, this
  assumption is exactly what is needed so that, whenever
  $\beta = a {\alpha_i} + b {\alpha_j} \in \Delta_+$ and
  $d {\alpha_i} < \beta$ for some $d \in \R$, then
  $d = 1$.  As a result, using the Coxeter relations, it follows
  inductively on depth that, if ${\alpha_i} < \beta$ for some
  $i \in I$, then if $\gamma < \beta$ and $\gamma \in \Delta_+$ is not
  simple, we also have ${\alpha_j} < \gamma$ for some $j \in I$.
  Thus, for all $\beta\in \Delta_+$, there exists $i \in I$ with
  ${\alpha_i} \leq \beta$.
\end{proof}

\begin{proof}[Proof of Theorem \ref{genmt}]
  It will be convenient to think of $m({v})$ as being allowed to be
  infinite (infinite if and only if the set appearing in the right
  hand side is empty). Similarly, call the number of moves required to
  reach a forbidden configuration ``infinite'' if and only if a
  forbidden configuration cannot be reached.  We clearly have
  $m({v}) \geq 0$, and Lemma \ref{curlem} implies that $m({v}) = 0$ if
  and only if ${v}$ is forbidden. Thus, using induction, the theorem
  may be restated as: if ${v}$ is not forbidden, then for any vertex
  $i$ with ${v}_i < 0$, we have $m(f_i {v}) \in \{m({v}), m({v})-1\}$,
  and there exists at least one such $i$ with
  $m(f_i {v}) = m({v}) - 1$.  Here, $\infty +c := \infty$ for any
  finite $c$.

Suppose that $\alpha \in \Delta_+$ and $j \in I$ are such that
$\alpha \cdot {v} < -1$ and ${v}_j < 0$. If we fire $j$, then the set
$\{\beta \in \Delta_+: \beta \cdot {v} < -1 \}$ changes by applying
$s_j$.  Hence, $m(f_j {v}) \in \{m({v})-1, m({v}), m({v})+1\}$. In
particular, $m(f_j {v}) \geq m({v})-1$.

Suppose that $\alpha \in \Delta_+$ is such that
$\alpha \cdot {v} < -1$ and $\dep(\alpha) -1 = m({v})$, and let
$i \in I$ be such that $s_i \alpha < \alpha$.  Then, if ${v}_i \geq
0$,
then $s_i \alpha \cdot {v} \leq \alpha \cdot {v} < -1$, which would
contradict the minimality of the depth of $\alpha$.  Thus, ${v}_i < 0$, and it
follows that $m(f_i {v}) = m({v}) - 1$. So, there exists $i$ such that
$m(f_i {v}) = m({v}) - 1$.

Next, suppose that ${v}_i < 0$ and $s_i \alpha > \alpha$. Then,
$\alpha \cdot f_i {v} \leq s_i \alpha \cdot f_i {v} < -1$. As a
result, we have $m(f_i {v}) \in \{m({v}), m({v}) -1\}$.  Thus, for any
$i \in I$ such that ${v}_i < 0$, we have
$m(f_i {v}) \in \{m({v}), m({v})-1\}$.
\end{proof}

\begin{proof}[Proof of Theorem \ref{onedirthm}]
  If $\alpha \cdot {v} < -1$, and $s_i \alpha > \alpha$, then
  ${v}_i < 0$ implies that $s_i \alpha \cdot f_i {v} < -1$ as well.
  As a result, although firing $i$ does not simply change
  $$Y_{{v}} := \{\beta \in \Delta_+: \beta \cdot {v} < -1 \, \,
  \text{and} \, \, \beta > {\alpha_i} \text{ for some $i$} \}$$
  by applying $s_i$, we still have $Y_{f_i {v}} \subseteq s_i
  Y_{{v}}$, which is all we need.
\end{proof}

\begin{rem}
  Note that, as a corollary of Lemma \ref{curlem}, we see that, for a
  general Coxeter group $W$, vertex $i \in I$, and matrix $C$, the set
  $\{j \in I \mid \exists b \in \R, b {\alpha_j} \in W
  {\alpha_i}\}$
  is the set of vertices $j$ connected to $i$ by a sequence of edges
  $i' \mapsto j'$ corresponding to odd integers $n_{i', j'}$.  It is
  clear that all such $j$ are in the set; conversely, if an edge
  corresponding to an even integer or $\infty$ is required to connect
  $i$ to $j$, then if $w {\alpha_i} = b {\alpha_j}$, then by
  modifying the elements of $C$ corresponding to the edges with even
  $n_{i'j'}$, we would be able to change the value $b$ such that
  $b {\alpha_j} \in W {\alpha_i}$. But this is impossible, since
  $b = 1$ whenever $c_{i'j'} = c_{j'i'}$ for all odd $n_{i'j'}$,
  and symmetrizing the latter values of $C$ would rescale $b$ by
  a fixed amount independent of the other values of $C$ (and
  independent of $b$ itself).
\end{rem}

\bibliographystyle{amsalpha}
\bibliography{references}

\end{document}